\documentclass{amsart}

\newtheorem{theorem}{Theorem}[section]
\newtheorem{lemma}{Lemma}[section]
\theoremstyle{remark}

\newtheorem{remark}{Remark}[section]
\newtheorem{corollary}[theorem]{Corollary}

\theoremstyle{definition}

\newtheorem{example}[theorem]{Example}

\numberwithin{equation}{section}

\begin{document}

\title[Growth of
difference quotients and logarithmic derivatives] {On the growth
of logarithmic differences, difference quotients and logarithmic
derivatives of meromorphic functions$^\dagger$}

\author{Yik-Man Chiang}
\address{Department of Mathematics, Hong Kong University of Science and
Technology, Clear Water Bay, Kowloon, Hong Kong,\ P. R. China}
\email{machiang@ust.hk}
\thanks{This research was supported in part by the Research Grants Council
of the Hong Kong Special Administrative Region, China
(HKUST6135/01P and 600806). The second author was also partially supported by the National Natural Science
Foundation of China (Grant No. 10501044) and by the HKUST PDF
Matching Fund.
}

\author{Shao-Ji Feng}
\address{Academy of Mathematics and Systems Science, Chinese Academy
of Sciences, Beijing, 100080, P. R. China}
\email{fsj@amss.ac.cn}
\thanks{${}^\dagger$ \textit{Many main results in this paper were presented in the ``Computational Methods and Function Theory" meeting held in Joensuu, Finalnd, $13^{\textrm{th}}-17^{\textrm{th}}$ June, 2005.}
}

\subjclass[2000]{Primary 30D30, 30D35, 39A05;
}

\date{11 Oct 2007 and in revised form, 25 November 2008, to appear \textit{Trans. Amer. Math. Soc.}}


\keywords{Difference operators, Poisson-Jensen formula,
Wiman-Valiron theory,  integrable difference equations}

\begin{abstract}
A crucial ingredient in the recent discovery by Ablowitz, Halburd,
Herbst and Korhonen \cite{AHH}, \cite {HK-2} that a connection
exists between discrete Painlev\'e equations and (finite order)
Nevanlinna theory is an estimate of the integrated average of
$\log^+|f(z+1)/f(z)|$ on $|z|=r$. We obtained essentially the same
estimate in our previous paper \cite{CF} independent of Halburd et
al \cite{HK-1}. We continue our study in this paper by
establishing complete asymptotic relations amongst the logarithmic
differences, difference quotients and logarithmic derivatives for
finite order meromorphic functions. In addition to the potential
applications of our new estimates in integrable systems, they are
also of independent interest. In particular, our findings show
that there are marked differences between the growth of
meromorphic functions with Nevanlinna order less than and greater
than one. We have established a ``difference" analogue of the
classical Wiman-Valiron type estimates for meromorphic functions
with order less than one, which allow us to prove that all entire
solutions of linear difference equations (with polynomial
coefficients) of order less than one must have positive rational
order of growth. We have also established that any entire solution
to a first order algebraic difference equation (with polynomial
coefficients) must have a positive order of growth, which is a
``difference" analogue of a classical result of P\'olya.
\end{abstract}

\maketitle







\section{Introduction}
We first set up some notation. Let $\eta$ be a fixed, non-zero
complex number, $\Delta f(z)= f(z+\eta)-f(z)$, and $\Delta^n
f(z)=\Delta(\Delta^{n-1}f(z))$ for each integer $n\ge 2$. In order
to simplify our notation, we shall use the same notation $\Delta$
for both a general $\eta$ and when $\eta=1$. The context will make
clear which quantity is under discussion. Equations written with
the above difference operators $\Delta^n f(z)$ are difference
equations. Let $E$ be a subset on the positive real axis. We
define the \textit{logarithmic measure} of $E$ to be
\begin{equation}
\label{E:defn-log-meas}
    \textrm{lm}(E)=\int_{E\cap (1,\infty)}\frac{dr}{r}.
\end{equation}
A set $E\in (1,\infty)$ is said to have finite logarithmic measure
if $\textrm{lm}(E)<\infty$. We adopt a more flexible interpretation on the Bachmann-Landau``big-$O$" notation \cite[p. 11]{WW1927}  so that for a complex function $f(z)$, $f(z)=O(\psi(r))$ is interpreted throughout this paper to mean that there is an $r_0>0$ such that $|f(z)/\psi(r)|<K$ holds for some  $K>0$ and for all $r=|z|>r_0$.
\medskip

Recently, there has been a renewed interest in difference and
$q$-difference
 equations in the complex plane $\textbf{C}$ (\cite{AHH}--
 \cite{BL}, \cite{CF}--\cite{GRP}, \cite{HK-1}--\cite{HK-5},
 \cite{Hay05}, \cite{HKLRT}--\cite{HLRY},
\cite{IY2004}, \cite{Ram}, \cite{Ruij}), and in particular,
Ablowitz, Halburd and Herbst \cite{AHH} proposed to use the Nevanlinna
order \cite{Hay64} as a detector of \textit{integrability} (i.e.,
solvability) of non-linear second order difference equations in
\textbf{C} (see \cite{AHH}, \cite{HK-2}, \cite{HK-4}, \cite{CR06};
see also \cite{RGH}--\cite{RGTT} and \cite[pp. 261--266]{GLS}).
Their theory is in close spirit with the classical
\textit{Painlev\'e ODE test} in which the solutions to ordinary
differential equations obtained from known integrable non-linear
PDEs via similarity reduction have particularly ``good" singularity
structure in \textbf{C}. That is, ``good" singularity structure of
solutions to ODEs in \textbf{C} can be regarded as a manifestation
of the integrability of certain non-linear PDES (see
\cite{AC1991}, pp. 98--100). There are also some works that focus
more on the function theoretic aspects of difference operators.
Halburd and Korhonen established a version of Nevanlinna theory
based on difference operators \cite{HK-3}, Bergweiler and Langley
\cite{BL} considered zeros of difference operators, and Ishizaki
and Yanagihara \cite{IY2004} developed a difference version of
Wiman-Valiron theory for entire functions of small growth. Halburd
and Korhonen \cite{HK-1} and the authors \cite{CF} studied the
growth of $f(z+\eta)/f(z)$ independently. The growth of
$f(z+\eta)/f(z)$ is crucial in the theory (\cite{HK-2}) of using
the Nevanlinna order as a detector of integrability of non-linear
second order difference equations. In particular, the authors
obtained, in \cite{CF}, that for any finite order meromorphic function
$f(z)$ of order $\sigma$, then for each $\varepsilon>0$,
\begin{equation}
\label{E:pointwise-est}
    \left|\frac{f(z+\eta)}{f(z)}\right|\le
    \exp(r^{\sigma-1+\varepsilon})
\end{equation}
holds for all $r$ outside a set of finite logarithmic measure. The
result was shown  to be best possible in a certain sense (see
\cite[e.g. 2.8]{CF}). The inequality (\ref{E:pointwise-est}) is in direct
analogy with the classical logarithmic derivative estimate by
Gundersen \cite{Gund88},
\begin{equation}
\label{E:pointwise-der}
    \left|\frac{f^\prime(z)}{f(z)}\right|\le
    |z|^{\sigma-1+\varepsilon}
\end{equation}
for all $|z|$ sufficiently large and outside a set of finite
logarithmic measure, which has countless applications (see e.g.
\cite{Laine}). In this paper we shall establish some estimates
that compare the growths of the logarithmic difference $\log
f(z+\eta)/f(z)$, the difference quotient $\Delta f/f$ and that of
$f^\prime/f$, and their applications to difference equations.
These results are extended to higher order differences and higher
order derivatives. In addition, we shall also exhibit examples of
meromorphic functions showing that our estimates are best possible when
interpreted in an appropriate sense.
\medskip

We shall establish, as our first main result, that for any finite order
meromorphic function of order $\sigma$, then for any given
$\varepsilon>0$,
\begin{equation*}
    2\pi i
    n_{z,\eta}+\log\frac{f(z+\eta)}{f(z)}=\eta\frac{f^\prime(z)}{f(z)}
    +O\big(r^{\beta+\varepsilon}\big),
\end{equation*}
or its equivalent form
\begin{equation}
\label{E:difference-log-der-1}
    \frac{f(z+\eta)}{f(z)}=e^{\eta\frac{f^\prime(z)}{f(z)}+
    O(r^{\beta+\varepsilon})},
\end{equation}
holds outside a set of $|z|=r$ of finite logarithmic measure,
$n_{z,\eta}$ is an integer that depends on both $z$ and $\eta$,
$\beta<\lambda-1$ when $\lambda<1$ and $\beta\le\lambda-1$ when
$\lambda\ge 1$,  where $\lambda$ is the
 exponent of convergence of the zeros and poles of $f$ (Theorem \ref{T:Fund-thm-1}).
 The above result holds for all finite order
meromorphic functions. In particular, it is easy to see that the inequality
(\ref{E:pointwise-est}) (\cite{CF}) mentioned above follows easily
from (\ref{E:difference-log-der-1}).
\medskip

Our second main result shows that if the order $\sigma$ is less
than one, then one can ``remove" the exponential function on the
right side of (\ref{E:difference-log-der-1}). More precisely, we
have
\begin{equation}
\label{E:delta-deriv-1} \frac{\Delta^kf(z)}{f(z)}=\eta^k
\frac{f^{(k)}(z)}{f(z)}+O(r^{(k+1)(\sigma-1)+\varepsilon}),
\end{equation}
which holds outside a set of $|z|=r$ of finite logarithmic measure.
 Equation (\ref{E:delta-deriv-1}) is a consequence of
the repeated applications of itself when $k=1$ and the following
estimate:
\begin{equation}
\label{E:delta-deriv-2}
    \left|\frac{\Delta^k f(z)}{f(z)}\right|\le
    |z|^{k\sigma-k+\varepsilon},
\end{equation}
which again holds for any $\varepsilon>0$ and for $z$ outside a
set of finite logarithmic measure. We note that our
(\ref{E:delta-deriv-2}) is in direct analogy with Gundersen's
(\ref{E:pointwise-der}). Our method of proof depends heavily on
the Poisson-Jensen formula.

\medskip

If we assume that $f(z)$ is an entire function, then we can
establish, as a direct consequence of (\ref{E:delta-deriv-1}), a
\textit{difference} Wiman-Valiron estimate
\begin{equation}
\label{E:difference-WV-1}
    \frac{\Delta^kf(z)}{f(z)}=\eta^k \Big(\frac{\nu(r, f)}{z}\Big)^k+
    O(r^{k\sigma-k-\gamma+\varepsilon}),
\end{equation}
which holds again outside an exceptional set of finite logarithmic
measure, and where $\gamma$ is a small positive constant which will
be made clear later. The notation $\nu(r, f)$ in
(\ref{E:difference-WV-1}) denotes the usual central index (see
\S\ref{S:WV-results}, \cite{Val} and \cite{PS}) of $f(z)$.
Although the remainder of our (\ref{E:difference-WV-1}) differs from
that of the classical estimate that involves the derivatives instead
of differences of $f(z)$, it is sufficient for most of our
applications to difference equations in \S\ref{S:applications}.
\medskip

Ishizaki and Yanagihara \cite{IY2004} established a ``difference"
version of Wiman-Valiron theory by expanding the entire function not
in terms of the usual basis $\{z^n\}$, but in terms of
\textit{factorial series}, that is, series written in terms of the
basis $\{(z)_n\}$ (where $(z)_n=z(z+1)\cdots (z+k-1)$). In
particular, their results are stated in terms of the central index
$\nu^*(r, f)$ with respect to $\{(z)_n\}$ instead of the usual
central index $\nu(r, f)$. Thus they need to impose a strong
restriction that $f(z)$ has order strictly less than $1/2$. Since
our assumption of $f$ only requires $\sigma(f)<1$, so our results
have a greater flexibility. When the entire function is of order
larger than one, the relation (\ref{E:difference-WV-1}) can no
longer hold since one has the relation
(\ref{E:difference-log-der-1}) instead. Thus the order assumption of
our result is best possible. Our (\ref{E:difference-WV-1}) has an
added advantage that it only involves the usual central index
$\nu(r, f)$ instead of the $\nu^*(r, f)$ which could be more
difficult to calculate.
\medskip 

Bergweiler and Langley investigated the zero distribution of $\Delta^k f(z)$  for entire functions of order  $\sigma(f) <1$ in \cite{BL}.  They obtain the relation
\cite[Lemma 4.2]{BL}
\begin{equation}
\label{E:Langley-Bergweiler}
	\Delta^k f(z)\sim f^{(k)}(z), 
\end{equation}
outside some exceptional set.  Although in close scrutiny one could derive the error estimate as in our (\ref{E:delta-deriv-1}) from some of their estimates, no explicit error bound is given in their paper.\footnote{When J. K. Langley showed us their manuscript of \cite{BL} after our presentation during the CMFT meeting in June, 2005, as noted in a footnote on our page 1, we had already obtained the main results presented in this paper.} In fact, it is clear that both the objectives and methods of approach between Bergweiler and Langley's paper and ours are very different. In addition, we exhibit examples showing that the major estimates in this paper, including the equation (\ref{E:delta-deriv-1}), are the best possible in some senses (see Example \ref{T:1st-order-counter-eg}). Moreover, our ``logarithmic differences" approach is natural from the Nevanlinna value distribution theory viewpoint where the corresponding error bounds are also crucial in our applications to difference equations as will be discussed in the next paragraph.
\medskip

 Both of the estimates
(\ref{E:delta-deriv-2}) and (\ref{E:difference-WV-1}) allow us to
estimate the growth of solutions in \textbf{C} of linear difference
equations,
\begin{equation}
\label{E:discrete-Eq}
    P_n(z) \Delta^n f(z)+\cdots + P_1(z)\Delta f(z)+P_0(z)f(z)=0,
\end{equation}
where $P_j,\, j=0,\cdots, n$ are polynomials, as in the classical
case of linear differential equations with polynomial coefficients
(see \cite{Val}, Chap. IV). In particular, we show that if the orders
of growth of entire solutions to (\ref{E:discrete-Eq}) are strictly
less than one, then they are equal to a set of positive rational
numbers, which can be obtained from the gradients of the
 Newton-Puisseux diagram of the associated algebraic equations.
Besides, we shall also show that any entire solution $f$ to the
first order algebraic difference equation
\begin{equation}
\label{E:1st-ord-diff-eq}
    \Omega\big(z,\, f(z),\, \Delta f(z)\big)=0
\end{equation}
with polynomial coefficients must have a positive order of growth.
This result is in close analogy with a corresponding result of
P\'olya \cite{polya} when $\Delta f$ is replaced by $f'$ and
with the same conclusion.
\medskip

This paper is organized as follows. The main result on logarithmic
differences will be stated in \S\ref{S:log-results} and proved in
\S\ref{S:proof-Fund-thm-1}. The result for first order difference
and those of higher order differences will be stated and proved in
  \S\ref{S:1st-order-results} and
\S\ref{S:higer-order-results} respectively. In addition, examples
will be constructed in both of \S\ref{S:log-results} and
\S\ref{S:1st-order-results} showing that the corresponding main
results there are best possible in a certain sense. We shall
formulate difference Wiman-Valiron type estimates in
\S\ref{S:WV-results}. Applications of our main results to linear
and first order algebraic difference equations
 will be discussed in \S\ref{S:applications}. Finally,
we shall discuss different aspects of the main results obtained in
this paper in \S\ref{S:conclusion}.
\eject\vfill

\section{A Logarithmic Difference Result}\label{S:log-results}
\begin{theorem}
\label{T:Fund-thm-1} Let $f(z)$ be a meromorphic function of order
$\sigma=\sigma(f)<\infty$, and let $\lambda^\prime$ and
$\lambda^{\prime\prime}$ be, respectively, the exponent of
convergence of the zeros and poles of $f$. Then for any given
$\varepsilon>0$, there exists a set $E\subset (1,\infty)$ of
$|z|=r$ of finite logarithmic measure, so that
\begin{equation}
\label{E:Fund-thm-est-1}
    2\pi i
    n_{z,\eta}+\log\frac{f(z+\eta)}{f(z)}=\eta\frac{f^\prime(z)}{f(z)}
    +O\big(r^{\beta+\varepsilon}\big),
\end{equation}
or equivalently,
\begin{equation}
\label{E:Fund-thm-est-2}
    \frac{f(z+\eta)}{f(z)}=e^{\eta\frac{f^\prime(z)}{f(z)}+
    O(r^{\beta+\varepsilon})},
\end{equation}
holds for $r\not\in E\cup[0,1]$, where  $n_{z,\eta}$ in
\eqref{E:Fund-thm-est-1} is an integer depending on both
$z$ and $\eta$,
    $\beta=\max\{\sigma-2,\,2\lambda-2\}$ if $\lambda<1$ and
    $\beta=\max\{\sigma-2,\,\lambda-1\}$ if $\lambda\ge1$ and
    $\lambda=\max\{\lambda^\prime,\lambda^{\prime\prime}\}$.
\end{theorem}
\medskip

\begin{remark}
Since $2\sigma-2<\sigma-1$ when $\sigma<1$ and
$\max\{\sigma-2,\,\lambda-1\}\le \sigma-1$ when $\sigma\ge 1$, so
we can easily deduce the estimate (\ref{E:pointwise-est}) obtained
in \cite{CF} as a consequence of our Theorem
\ref{E:Fund-thm-est-1}.
\end{remark}
\medskip

We shall construct an example here showing that the remainder in
(\ref{E:Fund-thm-est-2}) is best possible when
$\lambda=\max\{\lambda^{\prime},\lambda^{\prime\prime}\}<1$.

\begin{example}
\label{T:Fund-thm-2}
    Let $\lambda<1$ be a positive number and define
\begin{equation}
\label{E:counter-eg-1}
    f(z)=e^{z^k}\prod_{n=1}^\infty\Big(1-\frac{z}{n^{1/\lambda}}\Big).
\end{equation}
If $\lambda=0$, then we choose $f(z)=e^{z^k}$. The function $f$
has order $k$ when $k\ge 1$, and it has order $\lambda$ when
$k=0$, and the exponent of convergence of the zeros of $f$ is
$\lambda$.

 Let
\begin{equation}
\label{E:exceptional-set-1}
    F=\bigcup_{n=1}^\infty
    \big[n^{\frac{1}{\lambda}}+\frac13 n^{\frac{1}{\lambda}-1},\,
    n^{\frac{1}{\lambda}}+\frac23 n^{\frac{1}{\lambda}-1}\big].
\end{equation}

It can be seen that $F$ has infinite logarithmic measure on the
positive real axis. Let us consider $x\in F$. Then there exists an
integer $m$ such that
\begin{equation}
\label{E:exceptional-set-dist-3}
    m^{\frac{1}{\lambda}}+\frac13 m^{\frac{1}{\lambda}-1}\le x \le m^{\frac{1}{\lambda}}+\frac23  m^{\frac{1}{\lambda}-1}.
\end{equation}
We have
\begin{equation}
-3m^{1-\frac{1}{\lambda}}\le \frac{1}{m^{\frac{1}{\lambda}}-x}\le
0,\end{equation}
\begin{equation}
0\le \frac{1}{(m+1)^{\frac{1}{\lambda}}-x}\le
3m^{1-\frac{1}{\lambda}},\end{equation} and
\begin{equation}
m\ge \Big(\frac{x}{2}\Big)^\lambda.
\end{equation}
One can find positive constants $\delta$ and $c$ so that
\begin{equation}
\label{E:log-3}
    \log (1+w)-w\ge cw^2
\end{equation}
holds for all $w\in [-\delta,\, \delta]$. We consider $x\in F$ large
enough such that $3m^{1-\frac{1}{\lambda}}<\delta$. Then it follows
from (2.6) and (2.7) that
\begin{equation}
\label{E:exceptional-set-dist-1}
    -\delta\le \frac{1}{n^{\frac{1}{\lambda}}-x}\le \delta
\end{equation}
holds for all $n\in\textbf{N}$. We now apply (2.9) and (2.8) to
get for $x$ large enough and in $F$,
\begin{equation}
\begin{split}
    \log \frac{f(x+1)}{f(x)}-\frac{f'(x)}{f(x)} &=
    (x+1)^k-x^k -kx^{k-1}\\
    &\qquad+ \sum_{n=1}^\infty
    \log\Big(1-\frac{1}{n^{\frac{1}{\lambda}}-x}\Big)+\frac{1}{n^{\frac{1}{\lambda}}-x}\\
    &\ge \frac{k(k-1)}{2}x^{k-2}\big(1+o(1)\big)+c\sum_{n=1}^\infty
    \frac{1}{(n^{\frac{1}{\lambda}}-x)^2}\\
    & \ge \frac{k(k-1)}{2}x^{k-2}\big(1+o(1)\big)+c
    \frac{1}{(m^{\frac{1}{\lambda}}-x)^2}\\
    & \ge  \frac{k(k-1)}{2}x^{k-2}\big(1+o(1)\big) +\frac94 c
    m^{2-\frac{2}{\lambda}}\\
    &\ge  \frac{k(k-1)}{2}x^{k-2}\big(1+o(1)\big) +\frac94 c 2^{2-2\lambda}
    x^{2\lambda -2}.
\end{split}
\end{equation}

\end{example}
\bigskip

    We can see from the above  example that  Theorem
\ref{T:Fund-thm-1} is best possible when $\lambda<1$, in the sense
that the exponents $\sigma-2=k-2$ and $2\lambda-2$ in 
(\ref{E:Fund-thm-est-1}) are attained, and hence they cannot be
improved.

\section{Proof of Theorem \ref{T:Fund-thm-1}}
\label{S:proof-Fund-thm-1}
\subsection{Preliminary results }


\begin{lemma}
\label{L:log-1}
    Let us define
\begin{equation}
    \log w=\log|w|+i\arg w, \quad\quad -\pi\leq \arg w<\pi
\end{equation}
to be the principal branch of the logarithmic function in the complex
plane. Then we have
\begin{equation}
\label{E:log-0}
    \log(1+w)=O(|w|)
\end{equation}
and
\begin{equation}
\label{E:log-1}
    \log(1+w)- w=O(|w|^{2}),
\end{equation}
for $|w+1|\geq 1$.
\end{lemma}

\begin{proof}
We shall omit their elementary proofs.
\end{proof}

\begin{remark}
\label{R:log-defn} We note that it is clear that 
(\ref{E:log-1}) remains valid when we define the principal branch
by $\log w=\log|w|+i\arg w$, $-\pi< \arg w \le \pi$ instead.
\end{remark}

We shall apply the above lemma to prove

\begin{lemma}
 \label{L:log-2}
  Let us assume that we choose the principal branch as in Lemma
\ref{L:log-1}. Then we have
\begin{equation}
\label{E:log-1.5}
    \log(1+w)=O(|w|)+O\Big(\Big|\frac{w}{1+w}\Big|\Big)
\end{equation}
and
\begin{equation}
\label{E:log-2}
    \log(1+w)-w=O(|w|^{2})+O\Big(\Big|\frac{w}{1+w}\Big|^{2}\Big)
\end{equation}
hold for all $w$ in $\bf C$.
\end{lemma}

\begin{proof}
    If $|w+1|\geq 1$, then (\ref{E:log-1.5}) easily follows from
(\ref{E:log-0}). If, however, $|w+1|<1$, then
\begin{equation}
\label{E:log-3}
    \Big|1+\frac{-w}{1+w}\Big|=\frac{1}{|1+w|}>1.
\end{equation}
Thus we have
\begin{equation}
\label{E:log-3.1}
    \log(1+w)=-\log\frac{1}{1+w}=O\Big(\Big|1-\frac{1}{1+w}\Big|\Big)
    =O\Big(\Big|\frac{w}{1+w}\Big|\Big).
\end{equation}
This proves (\ref{E:log-1.5}) for all $w\in\textbf{C}$.

 Similarly, if $|w+1|\ge 1$ we easily see that (\ref{E:log-2}) follows from
 (\ref{E:log-1}). On the other hand, if $|w+1|< 1$, then (\ref{E:log-3}) holds, and (\ref{E:log-1}) implies that
\begin{equation}
\label{E:proof-log-1-1}
\begin{split}
     \log(1+w)-\frac{w}{1+w}&=
     -\Big[\log\Big(1+\frac{-w}{1+w}\Big)-
     \Big(\frac{-w}{1+w}\Big)\Big]\\
     &=O\Big(\Big|\frac{w}{1+w}\Big|^{2}\Big).
\end{split}
\end{equation}
We note that we have used a slightly different definition of the
principal logarithm as mentioned in Remark \ref{R:log-defn} in
order to handle the right hand sides of (\ref{E:log-3.1}) and
(\ref{E:proof-log-1-1}) above.

 On the
other hand, we have
\begin{equation}
\label{E:proof-log-1-2}
\begin{split}
    \Big(\frac{w}{1+w}\Big)-w &=\frac{-w^{2}}{(1+w)}=
    -w\cdot \frac{w}{1+w}\\
    &=O(|w|^{2})+O\Big(\Big|\frac{w}{1+w}\Big|^{2}\Big).
\end{split}
\end{equation}
 Thus we obtain, from combining (\ref{E:proof-log-1-1}) and
 (\ref{E:proof-log-1-2}),
\begin{equation*}
\begin{split}
    \log(1+w)-w &= \log(1+w)+\Big[-\Big(\frac{w}{1+w}\Big)
    +\Big(\frac{w}{1+w}\Big)-w\Big]\\
    &=O(|w|^{2})+O\Big(\Big|\frac{w}{1+w}\Big|^{2}\Big)
\end{split}
\end{equation*}
for all $w\in\textbf{C}$.
\end{proof}

\medskip

\begin{lemma}[\cite{Nev70}, p. 163]
\label{L:Nevanlinna}
 Let $f(z)$ be a meromorphic function in the complex plane, not identically zero.
Let $(a_{\nu})_{\nu\in N}$ and $(b_{\mu})_{\mu\in N}$, be the
sequence of zeros and poles, with due account of multiplicity, of
$f(z)$ respectively. Then for $|z|<R<\infty$,
\begin{equation}
\label{E:Poisson-Jensen}
\begin{split}
    \log f(z)=\frac{1}{2\pi}\int_0^{2\pi}& \log |f(Re^{i\phi})|\,
    \frac{Re^{i\phi}+z}{Re^{i\phi}-z}\,d\phi
    -\sum_{|a_{\nu}|<R}\log\frac{R^2-\bar{a}_{\nu}z}{R(z-a_{\nu})}\\
    &+\sum_{|b_{\mu}|<R}\log\frac{R^2-\bar{b}_{\mu}z}{R(z-b_{\mu})}+iC,
\end{split}
\end{equation}
where
\begin{equation}
\label{E:log-constant-1}
     C=\arg f(0)-\sum_{|b_{\mu}|<R}\arg\Big(-\frac{R}{b_{\mu}}\Big)
     +\sum_{|a_{\nu}|<R}\arg\Big(-\frac{R}{a_{\nu}}\Big)+2m_z\pi.
\end{equation}
\end{lemma}

\begin{remark}
\label{R:Nevanlinna}
 We recover the classical Poisson-Jensen
formula by taking the real parts on both sides of
(\ref{E:Poisson-Jensen}). Note that $m_z\in \textbf{N}$ in
(\ref{E:log-constant-1}) depends on the choice of branch of the
logarithm functions of both  sides of (\ref{E:Poisson-Jensen}),
and so it may depend on $z$ (but being piecewise continuous).
\end{remark}

We also require the following classical Cartan lemma.

\begin{lemma}[\cite{Car}; see also \cite{Lev}]
\label{L:Cartan} Let $z_1, z_2, \cdots, z_p$ be any finite
collection of complex numbers, and let $B>0$ be any given positive
number. Then there exists a finite collection of closed disks $D_1,
D_2, \cdots, D_q$ with corresponding radii $r_1, r_2, \cdots, r_q$
that satisfy
\begin{equation*}
    r_1+r_2+\cdots+r_q=2B,
\end{equation*}
such that if $z\notin D_j$ for $j=1,2,\cdots,q$, then there is a
permutation of the points $z_1,z_2,\cdots,z_p$, say,
$\hat{z}_1,\hat{z}_2,\cdots,\hat{z}_p$, that satisfies
\begin{equation*}
    |z-\hat{z}_l|>B\frac{l}{p}, \quad\quad l=1, 2, \cdots, p,
\end{equation*}
where the permutation may depend on $z$.
\end{lemma}

\subsection{Proof of Theorem \ref{T:Fund-thm-1}}

\begin{proof} For $\eta\not=0$, let $R>|z|+|\eta|$. We have by
(\ref{E:Poisson-Jensen}) and (\ref{E:log-constant-1}) that
\begin{equation}
\label{E:Nevanlinna-2}
\begin{split}
    \log f(z+\eta)&=\frac{1}{2\pi}\int_0^{2\pi} \log |f(Re^{i\phi})|\,
    \frac{Re^{i\phi}+z+\eta}{Re^{i\phi}-z-\eta}\,d\phi
    -\sum_{|a_{\nu}|<R}\log\frac{R^2-\bar{a}_{\nu}(z+\eta)}
    {R(z+\eta-a_{\nu})}\\
    &\quad +\sum_{|b_{\mu}|<R}\log\frac{R^2-\bar{b}_{\mu}(z+\eta)}
    {R(z+\eta-b_{\mu})}+i C_{\eta},
\end{split}
\end{equation}
where
\begin{equation}
\label{E:log-constant-2}
     C_\eta=\arg f(0)-\sum_{|b_{\mu}|<R}\arg\Big(-\frac{R}{b_{\mu}}\Big)
     +\sum_{|a_{\nu}|<R}\arg\Big(-\frac{R}{a_{\nu}}\Big)+2m_{z,\,\eta}\pi.
\end{equation}
We subtract (\ref{E:Poisson-Jensen}) from
(\ref{E:Nevanlinna-2}) and keep in mind that the integers
$m_{z,\eta},\,m_z$ in each step of the calculations below depend on the
choice of the logarithms as mentioned in Remark \ref{R:Nevanlinna}.
This yields
\begin{equation}
\label{E:log-difference}
\begin{split}
    \log\frac{f(z+\eta)}{f(z)}&=\frac{1}{2\pi}\int_0^{2\pi}\log
    |f(Re^{i\phi})|\,\frac{2\eta Re^{i\phi}}
    {(Re^{i\phi}-z-\eta)(Re^{i\phi}-z)}\,d\phi\\
    &\quad +\sum_{|a_{\nu}|<R} \bigg[\log\frac{R^2-\bar{a}_{\nu}z}{R(z-a_{\nu})}
    -\log\frac{R^2-\bar{a}_{\nu}(z+\eta)} {R(z+\eta-a_{\nu})}\bigg]\\
    &\quad +\sum_{|b_{\mu}|<R}\bigg[\log\frac{R^2-\bar{b}_{\mu}(z+\eta)}
    {R(z+\eta-b_{\mu})}-\log\frac{R^2-\bar{b}_{\mu}z}{R(z-b_{\mu})}\bigg]
    +2\pi(m_{z,\eta}-m_z)i \\
    &=\frac{1}{2\pi}\int_0^{2\pi}\log
    |f(Re^{i\phi})|\,\frac{2\eta Re^{i\phi}}
    {(Re^{i\phi}-z-\eta)(Re^{i\phi}-z)}\,d\phi\\
    &\quad-\sum_{|a_{\nu}|<R}\log\bigg[\Big(\frac{R^2-\bar{a}_{\nu}(z+\eta)}
    {R^2-\bar{a}_{\nu}z}\Big)\Big(\frac{z-a_\nu}{z+\eta-a_\nu}\Big)\bigg]\\
    &\quad+ \sum_{|b_{\mu}|<R}\log\bigg[\Big(\frac{R^2-\bar{b}_{\mu}(z+\eta)}
    {R^2-\bar{b}_{\mu}z}\Big)\Big(\frac{z-b_\mu}{z+\eta-b_\mu}\Big)\bigg]
    +2\pi (m_{z,\eta}-m_z)i\\
    &=\frac{1}{2\pi}\int_0^{2\pi}\log
    |f(Re^{i\phi})|\,\frac{2\eta Re^{i\phi}}
    {(Re^{i\phi}-z-\eta)(Re^{i\phi}-z)}\,d\phi\\
    &\quad-\sum_{|a_{\nu}|<R}\log\Big(1-\frac{\bar{a}_{\nu}\eta}
    {R^2-\bar{a}_{\nu}z}\Big)+\sum_{|a_{\nu}|<R}
    \log \Big(1+\frac{\eta}{z-a_\nu}\Big)\\
    &\quad+ \sum_{|b_{\mu}|<R}\log\Big(1-\frac{\bar{b}_{\mu}\eta}
    {R^2-\bar{b}_{\mu}z}\Big)- \sum_{|b_{\mu}|<R}\log
    \Big(1+\frac{\eta}{z-b_\mu}\Big)+2\pi(m_{z,\eta}-m_z)i.
\end{split}
\end{equation}

Differentiating the Poisson-Jensen formula (\ref{E:Poisson-Jensen})
yields
\begin{equation}
\begin{split}
\label{E:differeniated-Pois-Jen}
    \frac{f^\prime(z)}{f(z)}&=\frac{1}{2\pi}\int_0^{2\pi}\log
    |f(Re^{i\phi})|\frac{2Re^{i\phi}}{ (Re^{i\phi}-z)^2}d\phi\\
    &\quad\quad+\sum_{|a_{\nu}|<R}\frac{\bar{a}_{\nu}}{R^2-\bar{a}_{\nu}z}-
    \sum_{|b_{\mu}|<R}\frac{\bar{b}_{\mu}}{R^2-\bar{b}_{\mu}z}\\
    &\quad\quad+\sum_{|a_{\nu}|<R}\frac{1}{z-a_{\nu}}-
    \sum_{|b_{\mu}|<R}\frac{1}{z-b_{\mu}}.
\end{split}
\end{equation}
We multiply the left hand side of (\ref{E:differeniated-Pois-Jen})
by $\eta$ and subtract the product from the left hand side of
(\ref{E:log-difference}) with the assumption $|z|=r<R-|\eta|$. This
gives
\begin{equation}
\label{E:proof-thm-2.1-1}
\begin{split}
    \log \frac{f(z+\eta)}{f(z)}-\eta\frac{f^\prime(z)}{f(z)} &=
    \frac{\eta^2}{\pi}\int_0^{2\pi}\log|f(Re^{i\phi})|\frac{Re^{i\phi}}{
    (Re^{i\phi}-z-\eta)(Re^{i\phi}-z)^2}\, d\phi\\
    &\quad\quad -\sum_{|a_{\nu}|<R}\Bigg[\log\Big(1-\frac{\bar{a}_{\nu}\eta}{R^2-\bar{a}_{\nu}z}\Big)
        +\frac{\bar{a}_{\nu}\eta}{R^2-\bar{a}_{\nu}z}\Bigg]\\
    &\quad\quad+ \sum_{|b_{\mu}|<R}\Bigg[\log\Big(1-\frac{\bar{b}_{\mu}\eta}{R^2-\bar{b}_{\mu}z}\Big)
        +\frac{\bar{b}_{\mu}\eta}{R^2-\bar{b}_{\mu}z}\Bigg]\\
    &\quad\quad+\sum_{|a_{\nu}|<R}\Bigg[\log\Big(1+\frac{\eta}{z-a_{\nu}}\Big)
        -\frac{\eta}{z-a_{\nu}}\Bigg]\\
    &\quad\quad- \sum_{|b_{\mu}|<R}\Bigg[\log\Big(1+\frac{\eta}{z-b_{\mu}}\Big)
        -\frac{\eta}{z-b_{\mu}}\Bigg]+2\pi (m_{z,\eta}-m_z)i.
\end{split}
\end{equation}

We distinguish two cases depending on the exponent of convergence
of the zeros and poles of $f(z)$ in the rest of the proof. We
first consider the case\hfill\break $\lambda=\max\{\lambda^\prime,
\lambda^{\prime\prime}\} <1$.

We may now choose the branches of the logarithms on the right hand
side of (\ref{E:proof-thm-2.1-1})  as in Lemma \ref{L:log-2}.
Therefore we obtain

\begin{equation}
\label{E:proof-thm-2.1-2}
\begin{split}
    \Big|2\pi(m_z-m_{z,\eta})i +\log \frac{f(z+\eta)}{f(z)} &-\eta\frac{f^\prime (z)}{f(z)}\Big|
    \leq \frac{R|\eta|^2}{\pi(R-r-|\eta|)^3}\int_0^{2\pi}\big|\log |f(Re^{i\phi})|\big|
    \,d\phi\\
    &+O\Bigg[\sum_{|a_{\nu}|<R}\Big(\Big|\frac{\bar{a}_{\nu}\eta}{R^2-\bar{a}_{\nu}z}\Big|^2
    +|\frac{\bar{a}_{\nu}\eta}{R^2-\bar{a}_{\nu}z-\bar{a}_{\nu}\eta}\Big|^2\Big)\\
    &+\sum_{|b_{\mu}|<R}\Big(\Big|\frac{\bar{b}_{\mu}\eta}{R^2-\bar{b}_{\mu}z}\Big|^2
    +|\frac{\bar{b}_{\mu}\eta}{R^2-\bar{b}_{\mu}z-\bar{b}_{\mu}\eta}\Big|^2\Big)\\
    & +\sum_{|a_{\nu}|<R}\Big(\Big|\frac{\eta}{z-a_{\nu}}\Big|^2+\Big|\frac{\eta}{z-a_{\nu}
    +\eta}\Big|^2\Big)\\
    &+\sum_{|b_{\mu}|<R}\Big(\Big|\frac{\eta}{z-b_{\mu}}\Big|^2
    +\Big|\frac{\eta}{|z-b_{\mu}+\eta}\Big|^2\Big)\Bigg].
\end{split}
\end{equation}
We first note that
\begin{equation}
\label{E:proof-thm-2.1-3}
    \int_0^{2\pi}\big|\log |f(Re^{i\phi})|\big|\, d\phi=
    2\pi\big(m(R,f)+m(R,{1}/{f})\big).
\end{equation}
We next estimate the remaining terms in (\ref{E:proof-thm-2.1-2}).
We have
\begin{equation}
\label{E:proof-thm-2.1-4}
\begin{split}
    \sum_{|a_{\nu}|<R}\Big(\Big|\frac{\bar{a}_{\nu}\eta}{R^2-\bar{a}_{\nu}z}\Big|^2
    &+\Big|\frac{\bar{a}_{\nu}\eta}{R^2-\bar{a}_{\nu}z-\bar{a}_{\nu}\eta}\Big|^2\Big)\\
    &\quad \leq \sum_{|a_{\nu}|<R}\frac{2|\eta|^2R^2}{(R^2-Rr-R|\eta|)^2}
    =\frac{2|\eta|^2}{(R-r-|\eta|)^2}\, n(R,{1}/{f})
\end{split}
\end{equation}
and
\begin{equation}
\label{E:proof-thm-2.1-5}
\begin{split}
    \sum_{|b_{\mu}|<R}\Big(\Big|\frac{\bar{b}_{\mu}\eta}{R^2-\bar{b}_{\mu}z}\Big|^2
    &+\Big|\frac{\bar{b}_{\mu}\eta}{R^2-\bar{b}_{\mu}z-\bar{b}_{\mu}\eta}\Big|^2\Big)\\
    &\leq \sum_{|b_{\mu}|<R}\frac{2|\eta|^2R^2}{(R^2-Rr-R|\eta|)^2}
    =\frac{2|\eta|^2}{(R-r-|\eta|)^2}\, n(R,f).
\end{split}
\end{equation}
On the other hand, it is elementary that when $R^{\prime}>R>1$, we
have
\begin{equation}
\label{E:proof-thm-2.1-6}
\begin{split}
    N(R^{\prime}, f) &\geq \int_R^{R^{\prime}}\frac{n(t,f)-n(0, f)}{t}dt+n(0,f)\log R^{\prime}\\
    &\geq n(R,f)\int_R^{R^{\prime}}\frac{dt}{t}-n(0,f)\int_R^{R^{\prime}}\frac{dt}{t}+n(0,f)\log
    R^{\prime}\\
    & \geq n(R,f)\frac{R^{\prime}-R}{R^{\prime}}.
\end{split}
\end{equation}
Then for $R>1$, we have
\begin{equation}
\label{E:proof-thm-2.1-7}
     n(R,f)\leq 2N(2R,f),
\end{equation}
 and similarly
\begin{equation}
\label{E:proof-thm-2.1-8}
    n(R,{1}/{f})\leq 2N(2R,{1}/{f}).
\end{equation}

We may now combine
(\ref{E:proof-thm-2.1-2})--(\ref{E:proof-thm-2.1-5}),
(\ref{E:proof-thm-2.1-7}) and (\ref{E:proof-thm-2.1-8}) to get
\begin{equation}
\label{E:proof-thm-2.1-9}
\begin{split}
    \Big|2\pi(m_z-m_{z,\eta})i+\log \frac{f(z+\eta)}{f(z)}-\eta\frac{f^\prime(z)}{f(z)}\Big|
    &=O\bigg[\frac{R}{(R-r-|\eta|)^3}\Big(m(R,f)+m(R,{1}/{f})\Big)\\
    &\quad +\frac{1}{(R-r-|\eta|)^2}\big(N(2R,f)+N(2R,{1}/{f})\big)\\
    &\quad +\sum_{|c_k|<R}\Big(\frac{1}{|z-c_k|^2}+\frac{1}{|z-c_k+\eta|^2}\Big)\bigg],
\end{split}
\end{equation}
where $(c_k)_{k\in N}=(a_{\nu})_{\nu\in N}\cup(b_{\mu})_{\mu\in
N}$. Let $r>\max\{|\eta|, \frac13\}$ and $R=3r$. Then it follows
from (\ref{E:proof-thm-2.1-9}) and the finite order $\sigma$
assumption on $f$ that
\begin{equation}
\label{E:proof-thm-2.1-10}
\begin{split}
    2\pi(m_z-m_{z,\eta})i &+ \log \frac{f(z+\eta)}{f(z)} -\eta\frac{f^\prime(z)}{f(z)}\\
    &=
    O(r^{\sigma-2+\varepsilon})+O\bigg(\sum_{|c_k|<3r}\frac{1}{|z-c_k|^2}
    +\frac{1}{|z-c_k+\eta|^2}\bigg)\\
    &= O(r^{\sigma-2+\varepsilon})+O\bigg(\sum_{|d_k|<4r}\frac{1}{|z-d_k|^2}\bigg),
\end{split}
\end{equation}
where $(d_k)_{k\in N}=(c_k)_{k\in N}\cup(c_k-\eta)_{k\in N}$.

We now estimate the second summand in (\ref{E:proof-thm-2.1-10}) via
the Cartan lemma in the spirit in \cite{Gund88}.  Let $n(t)$ denote
the number of the points $d_k$ that lie in $|z|<t$. Then
\begin{equation}
\label{E:proof-thm-2.1-11}
    n(t)\le n(t+|\eta|, f)+n\big(t+|\eta|,{1}/{f}\big).
\end{equation}
We suppose that $h$ is any fixed positive integer and that $z$ is
confined to the annulus
\begin{equation}
\label{E:proof-thm-2.1-12}
    4^h\leq |z|=r\leq 4^{h+1}.
\end{equation}
Set $p=n(4^{h+2})$, $B=\frac{4^h}{h^2}$, and apply Lemma
\ref{L:Cartan} to the points $d_1,d_2,\cdots, d_p$, to conclude that
there exists a finite collection of closed disks $D_1, D_2, \cdots,
D_q$, whose radii has a total sum equal to $2B$, such that if
$z\notin D_j$ for $j=1,2,\cdots,q$, then there is a permutation of
the points $d_1,d_2,\cdots,d_p$, say,
$\hat{d_1},\hat{d_2},\cdots,\hat{d_p}$, that satisfies
\begin{equation}
\label{E:proof-thm-2.1-13}
    |z-\hat{d_k}|>B\frac{k}{p}, \qquad k=1, 2, \cdots, p.
\end{equation}
We note here that $q$ and $D_1, D_2, \cdots, D_q$ depend on $p$ and
then depend on $h$. Hence if $z\notin D_j$ for $j=1,2,\cdots,q$, we
have from (\ref{E:proof-thm-2.1-12}) and (\ref{E:proof-thm-2.1-13})
that
\begin{equation}
\label{E:proof-thm-2.1-14}
\begin{split}
    \sum_{|d_k|<4r}\frac{1}{|z-d_k|^2} &\leq \sum_{k=1}^p\frac{1}{|z-d_k|^2}
    =\sum_{k=1}^p\frac{1}{|z-\hat{d_k}|^2}\\
    &\leq \sum_{k=1}^p \frac{p^2}{B^2k^2} \leq
    \frac{p^2}{B^2}\sum_{k=1}^{\infty}\frac{1}{k^2}\\
    &\leq\frac{16[n(16r)]^2\log^4r}{r^2 \log^4 4}\sum_{k=1}^{\infty}\frac{1}{k^2}.
\end{split}
\end{equation}
We deduce from (\ref{E:proof-thm-2.1-11}) and
(\ref{E:proof-thm-2.1-14}) that
\begin{equation}
\label{E:proof-thm-2.1-15}
    \sum_{|d_k|<4r}\frac{1}{|z-d_k|^2}=O(r^{2\lambda-2+\varepsilon})
    +O(r^{2\lambda'-2+\varepsilon}).
\end{equation}

It remains to consider the exceptional sets arising from the discs
in the Cartan lemma. For each $h$, we define (it has been mentioned
that $q$ and $D_1, D_2, \cdots, D_q$ depend on $h$)
\begin{equation}
\label{E:proof-thm-2.1-16}
     Y_h=\{ r:\ z\in \bigcup_{j=1}^qD_j\ \mbox{such that} \ |z|=r\}
\end{equation}
and
\begin{equation}
\label{E:proof-thm-2.1-17}
    E_h=Y_h \cap [4^h,4^{h+1}].
\end{equation}
Then
\begin{equation}
\label{E:proof-thm-2.1-18}
     \int_{E_h}1\,dx\leq\int_{Y_h}1\,dx\leq 4B=\frac{4^{h+1}}{h^2}.
\end{equation}
Set
\begin{equation}
\label{E:proof-thm-2.1-19}
    E=\big\{\bigcup_{h=0}^{\infty}E_h\big\} \cap(1,\infty).
\end{equation}
We deduce, by combining (\ref{E:proof-thm-2.1-10}) and
(\ref{E:proof-thm-2.1-15}), that
\begin{equation}
\label{E:proof-thm-2.1-20}
    2\pi(m_z-m_{z,\eta})i+\log\frac{f(z+\eta)}{f(z)}=\eta\frac{f'(z)}{f(z)}
    +O(r^{\beta+\varepsilon})
\end{equation}
holds for all $z$ satisfying $|z|=r\notin [0,1]\cup E$. The sizes of
the exceptional sets can be calculated from
(\ref{E:proof-thm-2.1-17})--(\ref{E:proof-thm-2.1-18}). We have
\begin{equation*}
    \int_E \frac{1}{t}\, dt\leq\sum_{h=1}^{\infty}\int_{E_h}\frac{1}{t}\,dt
    \leq 4\sum_{h=1}^{\infty}\frac{1}{h^2}<\infty;
\end{equation*}
that is, $E$ has finite logarithmic measure. This completes the
proof when $\lambda <1$.

We now consider the remaining case when
$\lambda=\max\{\lambda^\prime,\lambda^{\prime\prime}\}\ge 1$. We
shall appeal directly to (\ref{E:proof-thm-2.1-1}) with
$r=|z|<R-|\eta|$. This together with (\ref{E:log-1.5}) and
(\ref{E:proof-thm-2.1-3}) yield, as in the case above when
$\lambda<1$,
\begin{equation}
\label{E:proof-thm-2.1-21}
\begin{split}
    \Big|2\pi (m_z-m_{z,\eta})i+\log \frac{f(z+\eta)}{f(z)} &
    -\eta\frac{f^\prime(z)}{f(z)}\Big|\le
    O\bigg[\frac{R}{(R-r-|\eta|)^3}\Big(m(R,f)+m(R,{1}/{f})\Big)\bigg]\\
    &\quad+\sum_{|c_{\mu}|<R}\bigg[\Big|\log\Big(1-\frac{\bar{c}_{\mu}\eta}{R^2-\bar{c}_{\mu}z}\Big)\Big|
        +\Big|\frac{\bar{c}_{\mu}\eta}{R^2-\bar{c}_{\mu}z}\Big|\bigg]\\
    &\quad+
    \sum_{|c_{\mu}|<R}\bigg[\Big|\log\Big(1+\frac{\eta}{z-c_{\mu}}\Big)\Big|
        +\Big|\frac{\eta}{z-c_{\mu}}\Big|\Bigg]\\
    &\le O\bigg[\frac{R}{(R-r-|\eta|)^3}\Big(m(R,f)+m(R,{1}/{f})\Big)\bigg]\\
    &\quad+O\Bigg[\sum_{|c_{\mu}|<R}\bigg(\Big|\frac{\bar{c}_{\mu}\eta}{R^2-\bar{c}_{\mu}z}\Big|
        +\Big|\frac{\bar{c}_{\mu}\eta}{R^2-\bar{c}_{\mu}z-\bar{c}_{\mu}\eta}\Big|\bigg)\Bigg]\\
    &\quad+ O\Bigg[\sum_{|c_{\mu}|<R}\bigg(\Big|\frac{\eta}{z-c_{\mu}}\Big|
        +\Big|\frac{\eta}{z+\eta-c_{\mu}}\Big|\bigg)\Bigg].\\
\end{split}
\end{equation}
We then choose $R=3r$, $r>\max\{|\eta|,1/3\}$. We also choose
$|z|$ as in (\ref{E:proof-thm-2.1-12}), $p=n(4^{h+2})$. Hence
\begin{equation}
\label{E:proof-thm-2.1-22}
\begin{split}
    \sum_{|c_{\mu}|<R}\bigg(\Big|\frac{\bar{c}_{\mu}\eta}{R^2-\bar{c}_{\mu}z}\Big| &+
    \Big|\frac{\bar{c}_{\mu}\eta}{R^2-\bar{c}_{\mu}z-\bar{c}_{\mu}\eta}\Big|\bigg)
    \le
    \sum_{|d_{\mu}|<4r}
    \Big|\frac{\bar{d}_{\mu}\eta}{(4r)^2-\bar{d}_{\mu}z}\Big|\\
    &\le\sum_{\nu=1}^{p}\frac{|\eta|}{4r-r}=p\frac{|\eta|}{3r}
    =O(r^{\lambda-1+\varepsilon}).
\end{split}
\end{equation}
We apply the Cartan argument similar to 
(\ref{E:proof-thm-2.1-12})--(\ref{E:proof-thm-2.1-14}) as in the
previous case when $\lambda<1$ to obtain
\begin{equation}
\label{E:proof-thm-2.1-23}
\begin{split}
    \sum_{|c_{\mu}|<R}\bigg(\Big|\frac{\eta}{z-c_{\mu}}\Big|
        &+\Big|\frac{\eta}{z+\eta-c_{\mu}}\Big|\bigg)\le
    \sum_{|d_{\mu}|<4r}\Big|\frac{\eta}{z-d_{\mu}}\Big|\\
    &\le|\eta|\sum_{\nu=1}^{p}
    \frac{1}{|z-\hat{d}_k|}\le |\eta|\sum_{\nu=1}^{p}
    \frac{p}{B}\frac{1}{\nu}\\
    &\le
    |\eta|p\frac{h^2}{4^h}\sum_{\nu=1}^{p}\frac{1}{\nu}
    \le
    |\eta|\frac{4n(16r)}{r}\Big(\frac{\log r}{\log 4}\Big)^2\log n(16r)
    =O(r^{\lambda-1+\varepsilon}).
\end{split}
\end{equation}
Combining (\ref{E:proof-thm-2.1-21}), (\ref{E:proof-thm-2.1-22})
and  (\ref{E:proof-thm-2.1-23}) yields
\begin{equation}
\label{E:proof-thm-2.1-24}
\begin{split}
    \Big|2\pi (m_z-m_{z,\eta})i+\log \frac{f(z+\eta)}{f(z)}
    -\eta\frac{f^\prime(z)}{f(z)}\Big|&\le
    O(r^{\sigma-2+\varepsilon})
    +O\big(\sum_{|d_{\mu}|<4r}
    \Big|\frac{\bar{d}_{\mu}\eta}{R^2-\bar{d}_{\mu}z}\Big|\Big)\\
    &\quad +O\Big(\sum_{|d_{\mu}|<4r}\Big|\frac{\eta}{z-d_{\mu}}\Big|\Big)\\
    &\le
    O(r^{\sigma-2+\varepsilon})+O(r^{\lambda-1+\varepsilon}),
\end{split}
\end{equation}
which holds outside a set of finite logarithmic measure. This completes the proof for the case when $\lambda\ge 1$ and
hence that of the theorem.
\end{proof}

\section{First Order Difference Quotients Estimates}
\label{S:1st-order-results}

\subsection{Main Results}

\begin{theorem}
\label{T:1st-order-thm}
 Let $f$ be a meromorphic function of order $\sigma(f)=\sigma<1$, and let $\eta$ be a fixed,
 non-zero number. Then for any $\varepsilon>0$, there exists a set $E\subset (1,\infty)$ that
depends on $f$ and has finite logarithmic measure, such that for all
$z$ satisfying $|z|=r\notin E\cup [0,1]$, we have
\begin{equation}
\label{E:1st-order-est-1}
    \frac{\Delta f(z)}{f(z)}=\frac{f(z+\eta)-f(z)}{f(z)}
    =\eta\frac{f^\prime(z)}{f(z)}+O(r^{2\sigma-2+\varepsilon}).
\end{equation}
\end{theorem}
\medskip

\begin{remark} We note that when $\eta=1$, then (\ref{E:1st-order-est-1})
assumes the form
\begin{equation}
\label{E:1st-order-est-2}
    \frac{\Delta
    f(z)}{f(z)}=\frac{f(z+1)-f(z)}{f(z)}=\frac{f^\prime(z)}{f(z)}+O(r^{2\sigma-2+\varepsilon}).
\end{equation}

This estimate will be extended to higher order difference quotients
in the next section.
\end{remark}

\medskip

\medskip

We next consider an example $f(z)$ showing that when
$\sigma(f)=1/2$, then the exponent ``$2\sigma-2=-1$" that appears
in the error term of (\ref{E:1st-order-est-1}) is best possible.

\medskip
\begin{example}
\label{T:1st-order-counter-eg} Let
  \begin{equation}
\label{E:cosine}
    f(z)=\cos \sqrt{z},
\end{equation}
which is clearly an entire function of order $1/2$. Then there is
a set $F$ of positive real numbers of infinite logarithmic
measure, such that for all $x\in F$, we have
\begin{equation}
\label{E:1st-order-est-4}
    \Big|\frac{\Delta f(x)}{f(x)}-\frac{f^\prime(x)}{f(x)}\Big|\ge
    \frac{1}{16x}{\big(1+o(1)\big)}.
\end{equation}
\end{example}
\bigskip

We first compute the numerators on the left hand side of
(\ref{E:1st-order-est-4}) by using trigonometric identities. This
gives
\begin{equation}
\label{E:counter-eg-2-1}
\begin{split}
    \Delta f(z)-f^\prime(z) &
    =\cos\sqrt{z+1}-\cos\sqrt{z}+\frac{\sin \sqrt{z}}{2\sqrt{z}}\\
    &= -2\sin\Big(\frac{\sqrt{z+1}+\sqrt{z}}{2}\Big)\sin\Big(\frac{\sqrt{z+1}-\sqrt{z}}{2}\Big)
     +\frac{\sin \sqrt{z}}{2\sqrt{z}}\\
    &=\frac{\sin
    \sqrt{z}}{2\sqrt{z}}-\frac{1}{2\sqrt{z}}\sin\Big(\frac{\sqrt{z+1}+\sqrt{z}}{2}\Big)\\
    &\qquad +\frac{1}{2\sqrt{z}}\sin\Big(\frac{\sqrt{z+1}+\sqrt{z}}{2}\Big)
    -2\sin\Big(\frac{\sqrt{z+1}+\sqrt{z}}{2}\Big)
    \sin\Big(\frac{\sqrt{z+1}-\sqrt{z}}{2}\Big)\\
    &=\frac{1}{2\sqrt{z}}\cdot
    2\cos\Big(\frac{\sqrt{z}+\frac{\sqrt{z+1}+\sqrt{z}}{2}}{2}\Big)
    \sin\Big(\frac{\sqrt{z}-\frac{\sqrt{z+1}+\sqrt{z}}{2}}{2}\Big)\\
    &\qquad + \sin\Big(\frac{\sqrt{z+1}+\sqrt{z}}{2}\Big) \Big[\frac{1}{2\sqrt{z}}
    -2\sin\Big(\frac{\sqrt{z+1}-\sqrt{z}}{2}\Big)\Big]\\
    &=-\frac{1}{\sqrt{z}}\cos\Big(\frac{\sqrt{z+1}+3\sqrt{z}}{4}\Big)
        \sin\Big(\frac{1}{4\big(\sqrt{z+1}+\sqrt{z}\big)}\Big)\\
    &\qquad + \sin\Big(\frac{\sqrt{z+1}+\sqrt{z}}{2}\Big) \Big[\frac{1}{2\sqrt{z}}
    -2\sin\Big(\frac{1}{2(\sqrt{z+1}+\sqrt{z})}\Big)\Big].\\
\end{split}
\end{equation}
We now consider $z=x>0$ and let $x\to+\infty$. This yields
\begin{equation}
\label{E:counter-eg-2-1}
\begin{split}
    \Delta f(x)-f^\prime(x) &=-\frac{1}{\sqrt{x}}
    \cos\Big(\frac{\sqrt{x+1}+3\sqrt{x}}{4}\Big)
    \cdot \frac{\big(1+o(1)\big)}{8\sqrt{x}}\\
    &\quad+ \sin\Big(\frac{\sqrt{x+1}+\sqrt{x}}{2}\Big) \Big[\frac{1}{2\sqrt{x}}
    -\frac{1}{(\sqrt{x+1}+\sqrt{x})}+O\Big(\frac{1}{(\sqrt{x+1}+\sqrt{x})}\Big)^3\Big]\\
    &=-\frac{1}{8x}
    \cos\Big(\frac{\sqrt{x+1}+3\sqrt{x}}{4}\Big){\big(1+o(1)\big)}\\
    &\quad+\sin\Big(\frac{\sqrt{x+1}+\sqrt{x}}{2}\Big)
    \Big[\frac{1}{2\sqrt{x}(\sqrt{x+1}+\sqrt{x})^2}+O\Big(x^{-3/2}\Big)\Big]\\
    &=-\frac{1}{8x}{\big(1+o(1)\big)}
    \cos\Big(\frac{\sqrt{x+1}+3\sqrt{x}}{4}\Big)+O\Big(x^{-3/2}\Big).
\end{split}
\end{equation}
We now let the subset $F$ of \textbf{R} be in the form
\begin{equation}
\label{E:counter-eg-2-2}
    F=\bigcup_{n=1}^\infty \big[(2\pi n-\pi/3)^2,\, (2\pi
    n+\pi/3)^2],
\end{equation}
which clearly has infinite logarithmic measure. We notice that for
$x\in F$ the function $\cos \sqrt{x}$ satisfies the inequality
\begin{equation}
\label{E:counter-eg-2-3}
    \frac12\le \cos\sqrt{x}\le 1.
\end{equation}
Combining (4.6) and (\ref{E:counter-eg-2-3})
yields, when $x\in F$ and $x\to \infty$,
\begin{equation}
\label{E:counter-eg-2-4}
    \frac{\Delta f(x)}{f(x)}-\frac{f^\prime(x)}{f(x)}
    \le-\frac{1}{8x}{\big(1+o(1)\big)}
    \cos\Big(\frac{\sqrt{x+1}+3\sqrt{x}}{4}\Big)+O\Big(x^{-3/2}\Big).
\end{equation}
It remains to consider the leading term on the right hand side of
(\ref{E:counter-eg-2-4}). We further notice that
\begin{equation}
\label{E:counter-eg-2-5}
\begin{split}
    \cos\Big(\frac{\sqrt{x+1}+3\sqrt{x}}{4}\Big)-\cos\sqrt{x}&=
    -2\sin\Big(\frac{\frac{\sqrt{z+1}+3\sqrt{z}}{4}+\sqrt{z}}{2}\Big)
    \sin\Big(\frac{\frac{\sqrt{z+1}+3\sqrt{z}}{4}-\sqrt{z}}{2}\Big)\\
    &=O(1)\cdot \sin\Big(\frac{1}{8(\sqrt{x+1}+\sqrt{x})}\Big)\\
    &=o(1).
\end{split}
\end{equation}
Thus we deduce from (\ref{E:counter-eg-2-5}), for all $x\in F$, that
\begin{equation}
\label{E:counter-eg-2-6}
    \cos\Big(\frac{\sqrt{x+1}+3\sqrt{x}}{4}\Big)=
    \cos\sqrt{x}+o(1)\ge \frac12+o(1).
\end{equation}

We can easily see that (\ref{E:1st-order-est-4}) follows from
(\ref{E:counter-eg-2-6}) and (\ref{E:counter-eg-2-4}).
\bigskip

\subsection{Proof of Theorem \ref{T:1st-order-thm}}

We require the classical estimate of Gundersen mentioned in the
Introduction.

\begin{lemma}[\cite{Gund88}]
\label{L:Gundersen} Let $f$ be a meromorphic function of finite
order $\sigma(f)=\sigma$. Then for any $\varepsilon>0$, there exists
a set $E\subset (1,\infty)$ that depends on $f$ and has finite
logarithmic measure, such that for all $z$ satisfying $|z|=r\notin
E\cup [0,1]$, we have
\begin{equation}
\label{E:Gundersen}
    \left|\frac{f^\prime (z)}{f(z)}\right|\le |z|^{\sigma-1+\varepsilon}.
\end{equation}
\end{lemma}
\bigskip

\begin{proof}
Given an arbitrary $\varepsilon$ such that $0<\varepsilon<1-\sigma$,
Theorem \ref{T:Fund-thm-1} and Lemma \ref{L:Gundersen} imply that
there exists a set $F\subset (1,\infty)$ that depends only on $f$
and has finite logarithmic measure, such that for all $z$
satisfying $|z|=r\notin F\cup [0,1]$,
\begin{equation}
\label{E:Proof-thm-3.1-1}
    2 n_z \pi i+\log\frac{f(z+\eta)}{f(z)}=O(r^{\sigma-1+\varepsilon})=o(1).
\end{equation}
Then there exists a constant $A>0$, such that for all $z$ satisfying
$|z|=r\notin F\cup [0,A]$,
\begin{equation}
\label{E:Proof-thm-3.1-2}
    \Big|2 n_z \pi i+\log\frac{f(z+\eta)}{f(z)}\Big|\leq 1.
\end{equation}
    Since
\begin{equation}
\label{E:Proof-thm-3.1-3}
     e^w-1=w+O(w^2),\qquad |w|\leq 1,
\end{equation}
we therefore deduce from (\ref{E:Proof-thm-3.1-1}),
(\ref{E:Proof-thm-3.1-2}), (\ref{E:Proof-thm-3.1-3}) and
(\ref{E:Fund-thm-est-2}) that for all $z$ satisfying $|z|=r\notin
F\cup [0,A]$,
\begin{equation}
\label{E:Proof-thm-3.1-4}
\begin{split}
    \frac{f(z+\eta)-f(z)}{f(z)}&=e^{2n_z\pi i+\log\frac{f(z+\eta)}{f(z)}}-1\\
    &=2 n_z \pi i+\log\frac{f(z+\eta)}{f(z)}+O\Big(2 n_z \pi i
    +\log\frac{f(z+\eta)}{f(z)}\Big)^2\\
    &=\eta\frac{f^\prime (z)}{f(z)}+O(r^{2\sigma-2+\varepsilon})
    +O\big((r^{\sigma-1+\varepsilon})^2\big)\\
    &=\eta\frac{f'(z)}{f(z)}+O(r^{2\sigma-2+2\varepsilon}).
\end{split}
\end{equation}
Now let $E=F\cup [1, A]$. Then $E$ has finite logarithmic measure
and hence (\ref{E:1st-order-est-1}) holds for all $z$ satisfying
$|z|=r\notin E\cup [0,1]$.
\end{proof}

\bigskip

\section{Higher Order Difference Quotients
Estimates}\label{S:higer-order-results}

\subsection{Main Results}

\begin{theorem}
\label{T:higer-order-1}
  Let $f$ be a meromorphic function of order $\sigma(f)=\sigma<1$.  Then for each given
$\varepsilon>0$, and integers $0\le j<k$, there exists a set $E\subset (1,\infty)$
that depends on $f$, and it has finite logarithmic measure, such that
for all $z$ satisfying $|z|=r\notin E\cup [0,1]$, we have
\begin{equation}
\label{E:higer-order-2}
    \Big|\frac{\Delta^kf(z)}{\Delta^j
    f(z)}\Big|\le |z|^{(k-j)(\sigma-1)+\varepsilon}.
\end{equation}
\end{theorem}
\bigskip

We deduce the following consequence from (\ref{E:higer-order-2}).
This is a higher order extension of Theorem
\ref{E:1st-order-est-1}.

\begin{theorem}
\label{T:higer-order-2} Let $f$ be a meromorphic function of order
$\sigma(f)=\sigma<1$, and let $\eta$ be a fixed, non-zero number.
Then for any $\varepsilon>0$, $k\in N$, there exists a set
$E\subset (1,\infty)$ that depends on $f$ and has finite
logarithmic measure, such that for all $z$ satisfying $|z|=r\notin
E\cup [0,1]$, we have
\begin{equation}
\label{E:higer-order-3}
    \frac{\Delta^kf(z)}{f(z)}=\eta^k\frac{f^{(k)}(z)}{f(z)}
    +O(r^{(k+1)(\sigma-1)+\varepsilon}).
\end{equation}
\end{theorem}
\medskip

We need the following elementary lemma to prove the theorems.

\medskip

\begin{lemma}
\label{L:difference-order} Let $f$ be a finite order meromorphic
function. Then for each $k\in N$,
\begin{equation}
\label{E:difference-order}
    \sigma(\Delta^k f)\leq\sigma(f).
\end{equation}
\end{lemma}

\begin{proof} Let $\sigma=\sigma(f)$. It is sufficient to prove the case when $k=1$.
We recall that the authors prove
\begin{equation}
\label{E:chiang-feng}
    T\big(r,\, f(z+\eta)\big)\sim
    T\big(r,\,f(z)\big)+O(r^{\sigma-1+\varepsilon})+O(\log r)
\end{equation}
in (\cite[Thm. 1]{CF}). Thus, we see immediately from
(\ref{E:chiang-feng}) that $\sigma(f(z+\eta))=\sigma(f(z))$ holds.
Hence
\begin{equation}
     \sigma(\Delta f(z))=\sigma(f(z+1)-f(z))\leq
\max\{\sigma(f(z+1)), \sigma(f(z))\}= \sigma(f(z)).
\end{equation}
Repeating the above argument $k-1$ times would yield (\ref{E:difference-order}).
\end{proof}
\medskip

\begin{remark} It is easy to see that the inequality
$\sigma(f)\le\sigma(\Delta^k f)$ does not hold for general
meromorphic functions. For example, one can take
$g(z)=\Gamma^\prime(z)/\Gamma(z)$.
\end{remark}

\subsection{Proof of Theorem \ref{T:higer-order-1}}

\begin{proof} Let $j$ be a non-negative integer amongst the set $\{0,\,1,\,\cdots, k-1\}$. We deduce from (\ref{E:1st-order-est-1}), (\ref{E:Gundersen}) and (\ref{E:difference-order})
that for each such $j$, there exists a set
$E_j\subset (1,\infty)$  that depends on $f$ and has finite
logarithmic measure, such that for all $z$ satisfying $|z|=r\notin
E_j\cup [0,1]$,
\begin{equation}
\label{E:proof-thm-5.1-1}
    \Big|\frac{\Delta^{j+1}f(z)}{\Delta^jf(z)}\Big|
    \le |z|^{\sigma_j-1+\varepsilon}\le |z|^{\sigma-1+\varepsilon},
\end{equation}
where $\sigma_j=\sigma(\Delta^jf)$. Now let $E=\bigcup_{j=0}^{k-1}E_j$,
which clearly has finite logarithmic measure. Hence, we have, after applying (\ref{E:proof-thm-5.1-1}) repeatedly, that
\begin{equation}
\label{E:proof-thm-5.1-2}
    \Big|\frac{\Delta^kf(z)}{f(z)}\Big|=\Big|\frac{\Delta^k f(z)}{\Delta^{k-1}f(z)}
    \cdot
    \frac{\Delta^{k-1} f(z)}{\Delta^{k-2} f(z)}\cdots \frac{\Delta f(z)}{f(z)}
    \Big|\le \big(|z|^{\sigma-1+\varepsilon}\big)^k
    \le |z|^{k\sigma-k+k\varepsilon}
\end{equation}
holds for all $z$ satisfying $|z|=r\notin E\cup [0,1]$. 
 Since $\varepsilon$ is arbitrary, this completes the proof of
 (\ref{E:higer-order-2}) when $j=0$. Let $G=\Delta^{j} f\, (j < k)$. Then Lemma \ref{L:difference-order}
  implies that $\sigma(G)\le \sigma(f)$. Thus
  (\ref{E:proof-thm-5.1-2})
  asserts that there is an exceptional set $E$ of finite logarithmic
  measure such that for all $z$ satisfying $|z|=r\not\in E$, we have
\begin{equation*}
    \Big|\frac{\Delta^k f(z)}{\Delta^jf(z)}\Big|=
    \Big|\frac{\Delta^{k-j} G(z)}{G(z)}\Big|\le
    |z|^{(k-j)(\sigma-1)+\varepsilon}.
\end{equation*}
This completes the proof of the theorem.
\end{proof}
\bigskip

\subsection{Proof of Theorem \ref{T:higer-order-2}}
\begin{proof}
The case when $k=1$ is just equation (\ref{E:1st-order-est-1}) in
Theorem \ref{T:1st-order-thm}. We shall prove 
(\ref{E:higer-order-3}) by induction on $k$. We assume that the
(\ref{E:higer-order-3}) is true for all $k$, $k\leq j$. We deduce
from Theorem \ref{E:1st-order-est-1}, Lemma \ref{L:Gundersen},
Theorem \ref{T:higer-order-1}, Lemma \ref{L:difference-order}, and the
fact that $\sigma(f')=\sigma(f)$, that there exists a set
$E\subset (1,\infty)$ that depends on $f$, which has finite
logarithmic measure, such that for all $z$ satisfying $|z|=r\notin
E\cup [0,1]$,
\begin{equation*}
\begin{split}
    \frac{\Delta^{j+1}f(z)}{f(z)}&=
    \frac{\Delta^{j+1}f(z)}{\Delta^jf(z)}\cdot\frac{\Delta^jf(z)}{f(z)}
    =\Big[\frac{(\Delta^jf)'(z)}{\Delta^jf(z)}+O(r^{2\sigma(f)-2+\varepsilon})\Big]
    \cdot\frac{\Delta^jf(z)}{f(z)}\\
    &=\frac{\Delta^j(f'(z))}{\Delta^jf(z)}\cdot\frac{\Delta^jf(z)}{f(z)}
    +O(r^{2\sigma(f)-2+\varepsilon})\cdot\frac{\Delta^jf(z)}{f(z)}\\
    &=\frac{\Delta^j(f'(z))}{f(z)}+O(r^{(j+2)(\sigma(f)-1)+2\varepsilon})\\
    &=\frac{\Delta^j(f'(z))}{f'(z)}\cdot\frac{f'(z)}{f(z)}+O(r^{(j+2)(\sigma(f)-1)+2\varepsilon})\\
    &=\Big[\frac{(f')^{(j)}(z)}{f'(z)}+O(r^{(j+1)(\sigma(f^\prime)-1)+\varepsilon})\Big]
    \cdot\frac{f'(z)}{f(z)}
    +O(r^{(j+2)(\sigma(f)-1)+2\varepsilon})\\
    &=\frac{f^{(j+1)}(z)}{f'(z)}\cdot\frac{f'(z)}{f(z)}+O(r^{(j+1)(\sigma(f)-1)+\varepsilon})
    \cdot\frac{f'(z)}{f(z)}+O(r^{(j+2)(\sigma(f)-1)+2\varepsilon})\\
    &=\frac{f^{(j+1)}(z)}{f(z)}+O(r^{(j+2)(\sigma(f)-1)+2\varepsilon}).
\end{split}
\end{equation*}
That is, the case $k=j+1$ is true, and the proof is complete.
\end{proof}

\section{Difference Wiman-Valiron type
Estimates}\label{S:WV-results}

Let $f(z)=\sum_{n=0}^{\infty}a_nz^n$ be an entire function in the complex
plane. If $r>0$, we define the maximum modulus $M(r, f)$ and
\textit{maximal term} $\mu(r, f)$ of $f$ by
\begin{equation*}
     M(r, f)=\max_{|z|=r}|f(z)|\qquad\textrm{and}\qquad \mu(r, f)=\max_{n\geq
     0}|a_n|r^n,
\end{equation*}
respectively. The \textit{central index} $\nu(r, f)$ is the greatest
exponent $m$ such that
\begin{equation*}
    |a_m|r^m=\mu(r, f).
\end{equation*}
We note that $\nu(r, f)$ is a real, non-decreasing function of
$r$.

We have the following fundamental result (see \cite{Laine} and
\cite{PS}) that relates the finite order of $f$ and its central
index.
\begin{lemma}
\label{L:WM-lemma-1}  If $f\not\equiv \textrm{const.}$ is an
entire function of order $\sigma$, then
\begin{equation}
\label{E:central-index}
     \limsup_{r\rightarrow\infty}\frac{\log\nu(r, f)}{\log r}
     =\sigma.
\end{equation}
\end{lemma}

We next quote the classical result of Wiman-Valiron (see also
\cite{Hay73}) in the form

\begin{lemma}[He and Xiao {\cite[pp.
28--30]{HS1988}}]
\label{L:WM-lemma-2} Let $f$ be a transcendental entire function. Let
$0<\varepsilon<\frac18$ and $z$ be such that $|z|=r$ and that
\begin{equation}
\label{E:Wiman-Valiron-1}
     |f(z)|>M(r, f)\big(\nu(r, f)\big)^{-\frac18+\varepsilon}
\end{equation}
holds. Then there exists a set $E\subset (1,\infty)$ of finite
logarithmic measure, such that
\begin{equation}
\label{E:Wiman-Valiron-2}
    \frac{f^{(k)}(z)}{f(z)}=\Big(\frac{\nu(r,
    f)}{z}\Big)^k\big(1+R_k(z)\big),
\end{equation}
\begin{equation}
\label{E:Wiman-Valiron-3}
    R_k(z)=O\big((\nu(r,f))^{-\frac18+\varepsilon}\big)
\end{equation}
holds for all  $k\in \text{\bf N}$ and all $r\notin E\cup [0,1]$.
\end{lemma}

We deduce from Theorem \ref{T:higer-order-2},
(\ref{E:Wiman-Valiron-2}) and (\ref{E:Wiman-Valiron-3}) in Lemma
\ref{L:WM-lemma-2}, the following.

\begin{theorem}
\label{T:difference-WV} Let $f$ be a transcendental entire
function of order $\sigma(f)=\sigma<1$, let
$0<\varepsilon<\frac18$ and $z$ be such that $|z|=r$, where
\begin{equation}
\label{E:max-place}
    |f(z)|>M(r, f)(\nu(r, f))^{-\frac18+\varepsilon}
\end{equation}
holds. Then for each positive integer $k$, there exists a set
$E\subset (1,\infty)$ that has finite logarithmic measure, such
that for all $r\notin E\cup [0,1]$,
\begin{equation}
\label{E:WV-zero-order}
    \frac{\Delta^kf(z)}{f(z)}=\eta^k\Big(\frac{\nu(r, f)}{z}\Big)^k
    \big(1+O((\nu(r, f))^{-\frac18+\varepsilon})\big),\quad \textrm{if}\quad \sigma=0,
\end{equation}
\begin{equation}
\label{E:WV-pos-order}
    \frac{\Delta^kf(z)}{f(z)}=\eta^k\Big(\frac{\nu(r, f)}{z}\Big)^k
    +O(r^{k\sigma-k-\gamma+\varepsilon}), \quad \textrm{if}\quad 0<\sigma<1,
\end{equation}
where $\gamma=\min\{\frac18\sigma, 1-\sigma\}$.
\end{theorem}

\section{Applications to Difference
Equations}\label{S:applications}

\begin{theorem}
\label{T:application-0}
  Let $P_0(z),\,\cdots,P_n(z)$ be polynomials such that
\begin{equation}
\label{E:appl-1-assumption}
    \max_{1\le j\le n}\{\deg P_j\}\le \deg(P_0).
\end{equation}
Let $f(z)$ be a meromorphic solution to the difference equation
\begin{equation}
\label{E:difference-eq-1}
    P_n(z)\Delta^n f(z)+\cdots+P_{1}(z)\Delta f(z)+P_0(z)f(z)=0.
\end{equation}
Then $\sigma (f)\ge 1$.
\end{theorem}

\begin{proof} Let us first assume that  a meromorphic solution $f$
has order $\sigma(f)<1$. We now write equation
(\ref{E:difference-eq-1}) in the form
\begin{equation}
\label{E:difference-eq-1.1}
    -1=\frac{P_n(z)}{P_0(z)}\frac{\Delta^n f(z)}{ f(z)}+\cdots+
    \frac{P_{1}(z)}{P_0(z)}\frac{\Delta f(z)}{f(z)},
\end{equation}
and choose $\varepsilon$ such that $0<\varepsilon<1-\sigma(f)$.
Theorem \ref{T:higer-order-1} asserts that there exists a set
$E\subset (1,+\infty)$ of finite logarithmic measure, such that for
all $|z|=r\not\in E\cup [0,1]$,
\begin{equation}
\label{E:difference-eq-3}
    \frac{\Delta^j f(z)}{f(z)}=o(1),
\end{equation}
for $1\le j\le n$. But
\begin{equation}
\label{E:difference-eq-4}
    \frac{P_j(z)}{P_0(z)}=O(1),
\end{equation}
as $|z|\to \infty$, for $1\le j\le n$. Therefore, if we choose
$|z|=r\not\in E\cup[0,1]$ and $|z|\to \infty$, then it follows from
(\ref{E:difference-eq-3}) and (\ref{E:difference-eq-4}) that we
have a contradiction in (\ref{E:difference-eq-1.1}).
\end{proof}

We note that this generalizes our earlier result \cite[Thm. 9.4]{CF},
which we state as
\begin{corollary}\label{T:appl3}
Let $Q_0(z),\, \cdots Q_n(z)$ be polynomials such that there exists
an integer $\ell,\, 0\le \ell\le n$ so that
\begin{equation}
\label{E:appl3-1}
  \deg (Q_\ell)> \max_{\substack{0\le \ell\le n\\ j\not=\ell}}\{\deg(Q_j)\}
\end{equation}
holds. Suppose $f(z)$ is a meromorphic solution to
\begin{equation}
\label{E:appl3-higher-order-linear}
    Q_n(z)y(z+n)+\cdots+Q_{1}(z)y(z+1) +Q_0y(z)=0.
\end{equation}
Then we have $\sigma(f)\ge 1$.
\end{corollary}

This is because if there is a polynomial coefficient of highest
degree amongst the $Q_0(z),\, \cdots Q_n(z)$, then the relation
\begin{equation}
\label{E:eg3-2}
    y(z+L)=\sum_{j=0}^L\binom{L}{j} \Delta^j y(z),\qquad
    L=0,\cdots, n,
\end{equation}
implies that we can transform
(\ref{E:appl3-higher-order-linear}) to (\ref{E:difference-eq-1}).
Then the corresponding $P_0$ obtained in (\ref{E:difference-eq-1})  has the highest degree. That is, (\ref{E:appl-1-assumption}) holds. Therefore, the result of the
corollary follows from Theorem \ref{T:application-0}.\hfill\qed
\bigskip

It is well known that each entire solution of the linear differential
equation
\begin{equation}
\label{E:linear-derivative}
    P_n(z) f^{(n)}(z)+\cdots+P_{1}(z) f^\prime(z)+P_0(z)f(z)=0,
\end{equation}
with polynomial coefficients, has order of growth equal to a
rational number which can be determined from the gradients of the
corresponding Newton-Puisseux diagram \cite{Hil76}. This classical
result can be proved from the Wiman-Valiron theory (Lemmas
\ref{L:WM-lemma-1} and \ref{L:WM-lemma-2}) (see also \cite[pp.
106--111, Appendix A]{Val} and \cite{GSW98}). We shall establish a
corresponding result for the linear difference equation
(\ref{E:difference-eq-1}) but only for entire solutions with order
strictly less than one. Our method is based on our ``difference"
Wiman-Valiron Theorem \ref{T:difference-WV}. Although the main
idea of our argument of using (\ref{E:WV-pos-order}) follows
closely that of the classical one for linear differential
equations, there are some details that require further
justifications. The reason is due to the different forms of the
remainders between (\ref{E:WV-pos-order}) in our Theorem \ref{T:difference-WV} and
the classical one in (\ref{E:Wiman-Valiron-2}).

\begin{theorem}
\label{T:application-1} Let $a_0(z),\,\cdots, a_n(z)$ be polynomial
coefficients of the difference equation
\begin{equation}
\label{E:difference-eq-2}
    a_n(z)\Delta^n f(z)+\cdots+a_{1}(z)\Delta f(z)+a_0(z)f(z)=0.
\end{equation}
Let $f$ be a transcendental entire solution of \eqref{E:difference-eq-2} with
$\sigma (f)< 1$. Then we have $\sigma(f)=\sigma=\chi$ where $\chi$
is a rational number which can be determined from a gradient of the
corresponding Newton-Puisseux diagram equation
\eqref{E:difference-eq-2}. In particular, $\chi>~0$.
\end{theorem}

\begin{proof} Let $f(z)$ be a transcendental entire solution of
(\ref{E:difference-eq-2}) with order $0<\sigma<1$. Lemma
\ref{L:WM-lemma-1} asserts that
\begin{equation}
\label{E:central-index-2}
     \limsup_{r\rightarrow\infty}\frac{\log\nu(r, f)}{\log r}
     =\sigma.
\end{equation}
This implies that for each $\delta>0$, there exists a sequence
$r_n\to \infty$, such that $r_{n+1}>r_n^{1+\delta}$ and $\nu(r_n,
f)\ge r_n^{\sigma-\delta}$.

We define
\begin{equation}
\label{E:central-index-3}
    F=\bigcup_{n=1}^\infty \big[r_n, \,r_n^{1+\delta}\big],
\end{equation}
which is a union of non-intersecting intervals and clearly has
infinite logarithmic measure. We consider all those $r$ that lie in
\begin{equation}
\label{E:central-index-4}
    r_{n}\le r \le r_{n}^{1+\delta}.
\end{equation}
Then since $\nu(r,f)$ is a non-decreasing function of $r$, so we
have trivially
\begin{equation}
\label{E:central-index-5}
    \nu(r,f)\ge \nu(r_{n},f)\ge r_n^{\sigma-\delta}
    \ge r^{\frac{\sigma-\delta}{1+\delta}}.
\end{equation}
We now choose $\delta$ and $\varepsilon$ so small that the
inequality
\begin{equation}
\label{E:central-index-6}
    k\sigma-k-\gamma+\varepsilon<\Big(\frac{\sigma-\delta}{1+\delta}-1\Big)k
\end{equation}
holds, where $\gamma$ is the constant in (\ref{E:WV-pos-order}).
Hence
\begin{equation}
\label{E:central-index-7}
    r^{k\sigma-k-\gamma+\varepsilon}=o\Big(\big|\frac{\nu(r,
    f)}{z}\Big|^k\Big),
\end{equation}
for $r$ in $F$. Let $z$ be such that (\ref{E:max-place}) holds.
Thus equation (\ref{E:WV-pos-order}) of Theorem \ref{T:difference-WV}
asserts that there exists a set $E$ of finite logarithmic measure
such that for all $|z|=r\in F\backslash \{E\cup[0,1]\}$, we have
\begin{equation}
\label{E:central-index-8}
    \frac{\Delta^kf(z)}{f(z)}=\eta^k\Big(\frac{\nu(r,
    f)}{z}\Big)^k\big(1+o(1)\big).
\end{equation}
Substituting (\ref{E:central-index-8}) into (\ref{E:difference-eq-2})
yields
\begin{equation}
\label{E:central-index-9}
    \sum_{j=0}^n \eta^j a_{A_j}^{(j)} \nu^j(r,\, f) z^{A_j-j} f(z)
    \big(1+o(1)\big)=0,
\end{equation}
where  $A_j=\deg a_j(z)$, $a_j(z)=\sum_{k=0}^{A_j}a_k^{(j)}z^k$.

Valiron \cite[p. 108 and Appendix A]{Val} asserts that we have,
for $r\in F\backslash \{E\cup[0,1]\}$, two real functions $B(r)$
and $\ell(r)$, where  $\ell(r)$ can only take values equal to the
gradients of the corresponding Newton-Puisseux diagram for
equation (\ref{E:difference-eq-2}) (which is always positive), and
$B(r)$ also can only assume a finite number of positive values
such that
\begin{equation}
\label{E:central-index-10}
    \lim_{r\to\infty}\frac{\nu(r, f)}{B(r)r^{\ell(r)}}=1.
\end{equation}
In particular, since $\ell(r)$ can only take a finite number of
rational values and $F\backslash \{E\cup[0, 1]\}$ has infinite
logarithmic measure, so there exists a rational $\chi>0$ such that the set
\begin{equation}
\label{E:central-index-11}
    \big\{r: \ell (r)=\chi,\quad r\in F\backslash
    \{E\cup[0,1]\}\big\}
\end{equation}
is unbounded. We conclude from (\ref{E:central-index-5}),
(\ref{E:central-index-2}), (\ref{E:central-index-10}) and
(\ref{E:central-index-11}) that
\begin{equation*}
    \frac{\sigma-\delta}{1+\delta}\le\chi\le\sigma.
\end{equation*}
Hence $\sigma=\chi$ by letting $\delta\to 0$.

It only remains to consider the case when $\sigma=0$. Since the
remainder of our ``difference" Wiman-Valiron estimate
(\ref{E:WV-zero-order}) is the same as the classical one
(\ref{E:Wiman-Valiron-2}), so we simply repeat the classical
argument as in \cite{Val}. This also shows that $\sigma=0$ must
be a gradient from the Newton-Puisseux diagram for equation
(\ref{E:difference-eq-2}),  which is impossible.
\end{proof}

We now consider a general first order algebraic difference equation.
P\'olya \cite{polya} proved the following classical theorem.
\begin{theorem}
\label{T:Polya}
    Let $f(z)$ be an entire solution of the first order algebraic
differential equation
\begin{equation}
\label{E:Polya-eqn}
    \Omega\big(z, f(z), f^\prime(z)\big)=0,
\end{equation}
with polynomial coefficients. Then $\sigma(f) >0$.
\end{theorem}

The proof utilizes the method of contradiction by assuming that if
$f(z)$ has order zero, then the central index $\nu(r, f)$ satisfies
\begin{equation}
\label{E:Polya-1}
    \lim_{r\to \infty} \frac{\nu(r, f)^k}{r}=0
\end{equation}
for every positive integer $k$ and (\ref{E:Wiman-Valiron-2}). Since our (\ref{E:WV-zero-order})
for the difference operator resembles (\ref{E:Wiman-Valiron-2}),
so P\'olya's argument applies verbatim to first order difference
equations also with polynomial coefficients (see also
\cite[Chap. 11]{Laine}). Thus we obtain
\begin{theorem}
\label{T:Polya-difference}
    Let $f(z)$ be an entire solution of the first order algebraic
difference equation
\begin{equation*}
\label{E:Polya-eqn}
    \Omega\big(z, f(z), \Delta f(z)\big)=0,
\end{equation*}
with polynomial coefficients. Then $\sigma(f) >0$.
\end{theorem}

\section{Concluding Remarks}\label{S:conclusion}
 Our method of proof of the main theorems depends heavily on the
Poisson-Jensen formula and is closely related to Nevanlinna theory. However, our results may be better
understood via the following formal approach to the difference
quotients and logarithmic derivatives. Let us denote by $E_\eta$   the operator $E_\eta f(z)=f(z+\eta)$, and define $D$ as $Df(z)=f'(z)$. Then we may write $\Delta f(z)=(E_\eta-1)f(z)$.
That is, $\Delta=E_\eta-1$ in the operator sense. Let us suppose in
addition that we have the expansion
\begin{eqnarray}
\label{E:taylor-expan-1}
    E_\eta f(z)&=& f(z+\eta)=f(z)+\eta
    f^\prime(z)+\frac{\eta^2}{2!}f^{\prime\prime}(z)+ \cdots\cr
    &=& f(z)+\eta
    Df(z)+\frac{\eta^2}{2!}D^2f(z)+ \cdots\cr
    &=& \left(1+\eta
    D+\frac{\eta^2D^2}{2!}+\frac{\eta^3D^3}{3!}+\cdots\right)f(z)\cr
    &=& (e^{\eta D})f(z).
\end{eqnarray}
That is, we have formally $E_\eta=e^{\eta D}$ (see \cite[p.
33]{Milne}). Substituting $\Delta=E_\eta-1$ into it yields again
formally
\begin{eqnarray}
\label{E:taylor-expan-2}
    \Delta^k &=& (e^{\eta D}-1)^k\cr
    &=& \Big( \eta D+\frac{\eta^2D^2}{2!}+\frac{\eta^3 D^3}{3!}+\cdots\Big)^k\cr
    &=& \eta^k D^k+ \frac{k}{2!}\,\eta^{k+1} D^{k+1}+\cdots.
\end{eqnarray}
Thus if we ignore the terms after the leading term in the last
line in (\ref{E:taylor-expan-2}), which was vigorously justified
in our argument in \S\ref{S:proof-Fund-thm-1}, we obtain
\begin{equation}
\label{E:taylor-expan-3}
    \frac{\Delta^k f}{f}=\eta^k \frac{f^{(k)}(z)}{f(z)}+\textrm{(small
    remainders)},
\end{equation}
which is formally the same as (\ref{E:delta-deriv-1}).

 The basic relations (\ref{E:difference-log-der-1}) and (\ref{E:delta-deriv-1}) give
 rise to
 the different growth patterns of finite order meromorphic
 functions of order greater than one and of order strictly less than one respectively.
 While the relations in both categories deserve more detailed study, the relation
 (\ref{E:difference-log-der-1}) for functions with
 order exactly one appears to be the most difficult to handle. But this
 category of functions includes many important functions, such as many classical
 special functions from mathematical physics.

Methods for solving linear difference equations have been
investigated by mathematicians as far back as the end of the
nineteenth century. For example, Milne-Thomson \cite[Chap.
XIV]{Milne} discussed methods for solving linear difference
equations in terms of factorial series via the operators ``$\pi$"
and ``$\rho$" that would ``converge everywhere, or nowhere or
converge in a half-plane in the right". It appears that the order
estimates for the meromorphic solution of Theorem \ref{T:application-0}
and those for the entire solutions to the linear difference
equations and first order non-linear difference equations in
Theorem \ref{T:application-1} and Theorem \ref{T:Polya-difference}
respectively are new. The topics deserve further study.

\vskip1cm \noindent\textbf{Acknowledgements} The authors would like to thank the referee for his/her constructive suggestions and for pointing out some errors in our original manuscript. These comments greatly improved the readability of the paper.

The second author thanks the Hong Kong University of Science and Technology for its hospitality during his visit from August 2004 to March 2005 where the main results presented in this paper were obtained.

\bibliographystyle{amsplain}

\end{document}